\documentclass[12pt]{article}
\usepackage{amsmath, amsthm}
\usepackage[utf8]{inputenc} 		
\usepackage[T1]{fontenc}
\usepackage[english]{babel}		
\usepackage{amsmath,amssymb,amsthm,tikz}	
\usepackage{graphicx}			
\usepackage{braids}
\usepackage{mathrsfs}
\usepackage{tikz-cd}
\usepackage{multicol}
\usepackage{hyperref}[]
\usepackage{mathtools}
\usepackage[left=1.6cm,right=1.5cm,top=2cm,bottom=2cm]{geometry}
\usepackage [all]{xy}
\graphicspath{{./graphics/}}
\usepackage{epstopdf}
\usetikzlibrary{matrix,arrows,decorations.pathmorphing} 
\usepackage{boxedminipage}
\newtheorem{theo}{Theorem}[section]

\newtheorem{prop}[theo]{Proposition}

\newtheorem{example}[theo]{Example}

\newtheorem{remark}[theo]{Remark}

\usepackage{color}			
\usepackage{courier} 
\newtheoremstyle{theorem}
{10pt} 
{10pt} 
{\sl} 
{\parindent} 
{\bf} 
{. } 
{ } 
{} 
\theoremstyle{theorem}

\newtheoremstyle{defi}
{10pt} 
{10pt} 
{\rm} 
{\parindent} 
{\bf} 
{. } 
{ } 
{} 
\theoremstyle{defi}

\begin{document} 
\date{} 
\title{On matrix representations for generalized braid groups}

\author{Abdoulrahim Ibrahim,\\[6pt]
IMAG, Univ Montpellier, CNRS.\\
Montpellier - 34090, France\\
e-mail: abdoulrahimibrahim@gmail.com\\[6pt]
}

\maketitle

\begin{abstract}
In the present paper, we construct a variant of the Burau representation of two generalizations of the classical braid group. For the Gassner representation, we propose an iterative procedure to find and generalize the extension of this representation.

\medskip



\end{abstract}

\section{Introduction}

In \cite{luddea}, a method for "lifting" representations from the Artin braid group on $n+1$ strands $B_{n+1}$ to the Artin braid group on $n$ strands $B_n$ is given by Lüdde and Toppan. 
Based on this method, we construct a variant of the Burau representation of the virtual braid group on $n$ strands $VB_n$, which is the subject of Theorem \ref{dernièreproposition}. 
\newline

The extension of Gassner representation of the pure braid group on $n$ strands $P_n$ was constructed independently by Bardakov \cite{barda05} and by Rubinsztein \cite{rubinsztein}. In the present paper, we propose an iterative procedure to find and generalize this representation, taking inspiration from the works \cite{mirko,luddea,luddeb}. We proceed as follows: we start with the fact that we give a subgroup $B$ of the automorphism group $Aut(F)$ of a free group $F$ of finite rank. So we have a natural embedding $B \hookrightarrow Aut(F)$ such that there is a canonical way to form a semi-direct product $F\rtimes B.$ Let $\tilde{\mathfrak{p}}:F \rtimes B \twoheadrightarrow B$ be the canonical projection whose its kernel is $F.$ We write $\mathfrak{p}:\mathbb{C}[F \rtimes B] \twoheadrightarrow \mathbb{C}[B].$ The kernel $\overline{I}_F$ of $\mathfrak{p}$ is the \emph{relative augmentation ideal} of $F$ in $\mathbb{C}[F \rtimes B].$ Using the embedding $B\hookrightarrow Aut(F),$ it is possible to define an action of $B$ on $\overline{I}_F.$ There is a representation of $B$ in terms of matrices whose entries are given in $\mathbb{C}[F \rtimes B].$ In Section \ref{sec:generalized Gassner}, we apply this procedure to the pure welded braid to obtain a faithful matrix representation of this group in Theorem \ref{4.3.6}, generalizing the Gassner representation of the pure Artin braid group.


\section{Two generalizations of the classical braid group}
\label{sec:two generalizations}
In this section we recall the standard definitions in terms of generators and relations of virtual and welded braid groups on $n$ strands. Let $VB_n$ denote the virtual braid group on the $n$ strands. This group first mentioned in \cite{kauffman} is the group defined by the presentation with generators $\sigma_1,\dots,\sigma_{n-1},\tau_1,\dots,\tau_{n-1}$ and the following relations
\begin{itemize}
\item[(V1)] $\sigma_{i}\sigma_{i+1}\sigma_{i}=\sigma_{i+1}\sigma_{i}\sigma_{i+1}$ for  $i= 1,2,\dots,n-2$
\item[(V2)] $\sigma_{i}\sigma_{j} = \sigma_{j}\sigma_{i}$ for  $|i-j| \geq 2$
\item[(V3)] $\tau_i^2=1$ for $i=1,\dots,n-1$
\item[(V4)] $\tau_{i}\tau_{i+1}\tau_{i}=\tau_{i+1}\tau_{i}\tau_{i+1}$ for  $i= 1,2,\dots,n-2$
\item[(V5)] $\tau_{i}\tau_{j} = \tau_{j}\tau_{i}$ for  $|i-j| \geq 2$
\item[(V6)] $\sigma_i \tau_j =\tau_j \sigma_i$ for $|i-j| \geq 2$
\item[(V7)] $\sigma_{i} \tau_{i+1}\tau_{i} = \tau_{i+1}\tau_{i} \sigma_{i+1}$ for $i= 1 \dots n-2$
\label{Eq:virtualrelations}
\end{itemize}
It is well-known that the relations (V1)-(V2) are defining relations of the classical braid group on $n$ strands, $B_n$, where the elements $\sigma_1,\dots,\sigma_{n-1}$ generate the classical braid group $B_n$ in $VB_n$. It is also known that the relations (V3)-(V5) are defining relations of the symmetric group $S_n$ on $n$ letters, where the elements $\tau_1,\dots,\tau_{n-1}$ generate the symmetric group $S_n$ in $VB_n.$ Hence $VB_n$ can be presented by the free product of $B_n$ and $S_n$ modulo relations (V6)-(V7). The mirror of the relation (V7) (i.e. when read backwards) hold in $VB_n:$
\begin{eqnarray}
\tau_i \tau_{i+1} \sigma_{i}=\sigma_{i+1}\tau_{i}\tau_{i+1} \text{ $ $ } (i=1,\dots,n-2).
\end{eqnarray}
But however, the following relation do not hold in $VB_n$ (see \cite{goussarovpolyakviro}):
\begin{eqnarray}
\tau_{i}\sigma_{i+1}\sigma_{i}=\sigma_{i+1}\sigma_{i} \tau_{i+1} \text{ $ $ } (i=1,\dots,n-2)
\label{Eq:forbiddenrelation1}
\end{eqnarray}
The welded braid group on $n$ strands $WB_n$ is by definition the quotient of $VB_n$ by this relation, usually called the forbidden relation. Thus, $WB_n$ is the group obtained from $VB_n$ by introducing the forbidden relation \eqref{Eq:forbiddenrelation1} and, moreover, contains important elements of the form
\begin{eqnarray}
\xi_{i,j} &=& \tau_{i}\tau_{i+1}\cdots \tau_{j-2}\tau_{j-1}\sigma_{j-1} \tau_{j-2}\cdots\tau_{i+1} \tau_{i} \text{  $ $   } (i<j)\nonumber\\
\xi_{i,j} &=& \tau_{i-1}\tau_{i-2}\cdots\tau_{j-2}\tau_{j-1}\sigma_j\tau_j\tau_{j-1} \cdots \tau_{i-1} \text{  $ $   } (j<i)\nonumber\\
\end{eqnarray}
The elements $\xi_{i,j}$ ($1\leq i\neq j\leq n$) generate the welded pure braid group $PW_n$ also known as the group of basis-conjugating automorphisms in \cite{savushikina}, or the McCool group in \cite{berceanupapadima} which is defined as the kernel of the projection $WB_n\twoheadrightarrow S_n$ by sending $\sigma_{i}\mapsto \tau_{i}$ and $\tau_{i}\mapsto \tau_{i}.$ In \cite{mccool}, McCool proved that the following relations known as the \emph{McCool relations} determine a presentation of $PW_n$:
\begin{eqnarray}
\left[\xi_{k,j}, \; \xi_{s,t}\right] &=&1 \text{ if } \; \{i,j\} \cap \{s,t\}= \emptyset \nonumber,\\
\left[\xi_{i,j}, \; \xi_{k,j}\right ]&=&1 \text{ for } \; i,j, k \text{ distinct }\nonumber,\\
\left[\xi_{i,j}\cdot \xi_{kj}, \; \xi_{i,k}\right ]&=&1 \text{ for } \; i,j, k  \text{ distinct }.
\label{Eq:2.7}	
\end{eqnarray}
where $\left[a,b\right]:=a^{-1}b^{-1}ab$. In \cite{savushikina}, it is proved that the welded braid group $WB_n$ splits as a semidirect product, $WB_n=PW_n \rtimes S_n$ and there is an analogous statement for the virtual braid group \cite{bardakov04}. Moreover, $WB_n$ has a ''twin`` group denoted by $\mathcal{WB}_n$. This is the group corresponding to $WB_n$ as a set, and we have all the relations as in $WB_n$ except that the forbidden relation \eqref{Eq:forbiddenrelation1} which is replaced by 
\begin{eqnarray}
\tau_{i+1} \sigma_{i} \sigma_{i+1}=\tau_{i}\sigma_{i+1}\sigma_{i}.
\end{eqnarray}
Both are isomorphic and the isomorphism is given by the map $\mathfrak{V}_n:\mathcal{WB}_n \rightarrow WB_n$ that sends $\sigma_{i}\mapsto \sigma^{-1}_{i}$ and $\tau_{i}\mapsto \tau_{i}.$ 
\begin{remark}
The relation $\tau_{i+1} \sigma_{i} \sigma_{i+1}=\tau_{i}\sigma_{i+1}\sigma_{i}$ which holds in $\mathcal{WB}_n$ does not hold in $WB_n.$
\end{remark}
We conclude this section with a representation for the welded braid group by automorphisms of the free group of rank $n.$ Recall that the subgroup of $PW_{n+1}$ generated by the elements $\xi_{n+1,1},\dots, \xi_{n+1,n}$ is a free group of rank $n,$ which we denote by $F_n$ (see \cite{barda03}). For $1\leq i \leq n- 1,$ let $\rho_i:F_n \rightarrow F_n$ and $\vartheta_i:F_n \rightarrow F_n$ be the automorphisms defined respectively by
\begin{eqnarray}
\rho_i: \left\{
\begin{array}{ll}
\xi_{n+1,i} &\mapsto \xi_{n+1,i} \xi_{n+1,i+1} \xi_{n+1,i}^{-1} \\
\xi_{n+1,i+1}&\mapsto \xi_{n+1,i}\\
\xi_{n+1,j} &\mapsto \xi_{n+1,j} $ $ \mbox{ $ $ }  \forall j\neq \{i,i+1\}.
\end{array}
\right. \text{ and  } \vartheta_i: \left\{
\begin{array}{ll}
\xi_{n+1,i}&\mapsto  \xi_{n+1,i+1} \\
\xi_{n+1,i+1} &\mapsto \xi_{n+1,i}\\
\xi_{n+1,j} &\mapsto \xi_{n+1,j} $ $ \mbox{ $ $ }  \forall j\neq \{i,i+1\}.
\end{array}
\right.
\end{eqnarray}
One can easily show the following.
\begin{prop}
The mapping $\sigma_i\mapsto \rho_i$ and $\tau_{i}\mapsto \vartheta_i$, $1\leq i\leq n-1$, determines a representation $\psi:WB_n\rightarrow Aut(F_n)$. 
\end{prop}
The above representation $\psi:WB_n\rightarrow Aut(F_n)$  is commonly known as the Artin representation and it is faithful (for a proof see \cite{frimrourke}).

\section{Generalized Burau and Gassner representation}
Here we construct matrices representations that generalize the classical Burau and Gassner representation.
\subsection{The virtual braid group lift}
\label{sec:Un relèvement de représentation de tresses virtuelles}

The idea for this section is based on \cite{luddea}. Let $\mathbb{C}$ denote the field of complex numbers. Let us now take $\mathbb{C}[WB_{n+1}]$ to be the group algebra of $WB_{n+1}$ and consider an abstract set $\Delta^{(1)}_n:=\{\delta^{(1)}_1,\dots,\delta^{(1)}_n\}$. Then we formally construct a left free module on $\mathbb{C}[WB_{n+1}]$ from this set. The free $\mathbb{C}[WB_{n+1}]$-module to the left of rank $n$ thus obtained is denoted by $\mathcal{F}^{(1)}_n$ and it is generated by the elements $\delta^{(1)}_1,\dots,\delta^{(1)}_n.$ This module carries the following representation of the virtual braid group $VB_n$ defined as an action to the right of $VB_n$ on $\mathcal{F}^{(1)}_n:$ the generators $\left(\sigma_i\right)^{n-1}_{i=1}$ and $\left(\tau_i\right)^{n-1}_{i=1}$ in this representation act according to
\begin{eqnarray}
\delta^{(1)}_j.\sigma_{i} =\sigma_{i} &\times& \left(1-\alpha\cdot\xi_{n+1,i} \xi_{n+1,i+1} \; \xi^{-1}_{n+1,i}\right)\cdot \delta^{(1)}_{i}+ \alpha\cdot\xi_{n+1,i}\cdot \delta^{(1)}_{i+1}  \text{  $ $   if $ $ } j=i \nonumber\\
\sigma_i & \times & \delta^{(1)}_i  \text{  $ $   if $ $ } j=i+1\nonumber \\
\sigma_i& \times & \delta^{(1)}_j   \text{  $ $   if $ $ } j\neq \{i,i+1\} \nonumber\\
\nonumber\\
\delta^{(1)}_j.\tau_{i}=\tau_{i} &\times& \beta^{-1} \cdot \delta^{(1)}_{i+1} \text{  $ $   if $ $ } j=i \nonumber\\
\tau_i&\times& \beta\cdot \delta^{(1)}_i   \text{  $ $   if $ $ } j=i+1 \nonumber \\
\tau_i&\times& \delta^{(1)}_j \text{  $ $   if $ $ } j\neq \{i,i+1\} \nonumber\\
\label{Eq:7}
\end{eqnarray}
where $\alpha$ and $\beta$ are nonzero complex parameters (not necessarily distinct). The action on the right \eqref{Eq:7} can be written as follows:
\begin{eqnarray}
\delta_j.\sigma_{i}&=&\sum_{k=1}^{n} \mathcal{V}(\sigma_i)^k_j \; \delta_k, \text{  $$ $$   } (1\leq j \leq n)\nonumber \\
\nonumber\\
\delta_j.\tau_{i}&=& \sum_{k=1}^{n} \mathcal{V}(\tau_i)^k_j \; \delta_k,\text{  $$ $$   } (1\leq j \leq n)  
\label{Eq:8}
\end{eqnarray}
with
\begin{eqnarray}
\mathcal{V}(\sigma_i)^k_j:=\bordermatrix {
	&  & i & & & i+1 &  \cr
	& \sigma_i\cdot I_{i-1} &  & & &  \cr
	&  & \sigma_i\cdot(1- \alpha\cdot\xi_{n+1,i} \cdot \xi_{n+1,i+1} \cdot \xi^{-1}_{n+1,i}) & & & \sigma_i \cdot \alpha\cdot\xi_{n+1,i} &  \cr
	&  &  & & & &  \cr
	&  & \sigma_i & & &0 &  \cr
	& & & & & & \cr
	&  & & & & &\sigma_i \cdot I_{n-i-1} \cr 
} \nonumber\\
\label{Eq:9}
\end{eqnarray}
and
\begin{eqnarray}
\mathcal{V}(\tau_i)^k_j:= \bordermatrix {
	&  & i & & i+1 &  \cr
	& \tau_i\cdot I_{i-1} & & &  &  \cr
	&  & 0 & &\beta^{-1}\cdot\tau_i &  \cr
	&  &  & & &  \cr
	&  & \beta\cdot\tau_i & &0&  &  \cr
	&  && & &\tau_i\cdot I_{n-i-1} \cr
} 
\label{Eq:10}.
\end{eqnarray}
where $\alpha,\beta \in \mathbb{C}^{\ast}.$ 
\begin{prop}
Let $\alpha,\beta$ be nonzero complex numbers with $\alpha \neq \beta.$ Then the following assignment defines a matrix representation of $VB_n$:
	\begin{eqnarray}
	\phi_n:VB_n  \longrightarrow  GL(n,\mathbb{C}[WB_{n+1}] )
	\end{eqnarray}	defined by
	\begin{eqnarray}
	\sigma_i \mapsto \mathcal{V}(\sigma_i):=\bordermatrix {
		&  & i & & & i+1 &  \cr
		& \sigma_i\cdot I_{i-1} &  & & &  \cr
		&  & \sigma_i\cdot (1- \alpha\cdot\xi_{n+1,i} \; \xi_{n+1,i+1} \; \xi^{-1}_{n+1,i}) & & & \alpha\cdot\sigma_i\cdot\xi_{n+1,i} &  \cr
		&  &  & & & &  \cr
		&  & \sigma_i & & &0 &  \cr
		& & & & & & \cr
		&  & & & & &\sigma_i\cdot I_{n-i-1} \cr 
	}  
	\end{eqnarray}
	and 
	\begin{eqnarray}
	\tau_{i} \mapsto \mathcal{V}(\tau_i):=\bordermatrix {
		&  & i & & i+1 &  \cr
		& \tau_i\cdot I_{i-1} & & &  &  \cr
		&  & 0 & &\beta^{-1}\cdot\tau_i &  \cr
		&  &  & & &  \cr
		&  &\beta\cdot \tau_i & &0&  &  \cr
		&  && & &\tau_i\cdot I_{n-i-1} \cr
	} .
	\end{eqnarray}
Moreover, the representation $\phi_n$ factors though $\mathcal{WB}_n$ if and only if $\beta=\alpha$ and $WB_n$ if and only if $\beta=1.$ 
	\label{theofaithful}
\end{prop}
We derive from Proposition \ref{theofaithful} a variant of the Burau representation of the virtual braid group, as far as I know, which seems to be new which we state in the following theorem.
\begin{theo}
Let $\alpha,\beta$ denote two nonzero complex parameters with $\alpha \neq \beta$ and $\beta \neq 1.$	The mapping $\tilde {h}:WB_{n+1}\rightarrow \mathbb{C}$ such that $\sigma_i \mapsto 1$ and $\tau_i \mapsto 1$, applied to each entry of the above matrix yields the variant of the Burau representation ${}^{\tilde{h}}\phi_n: VB_n \longrightarrow GL(n,\mathbb{C})$ defined by	
	\begin{eqnarray}
	\sigma_i  \mapsto	{}^ {\tilde{h}}  \mathcal{V}(\sigma_i)=\bordermatrix {
		&  & i & &i+1 &  \cr
		& I_{i-1} &&  &  &  \cr
		&  & 1-\alpha & &\alpha & \cr
		& & & & & \cr 
		&  & 1 & &0 &  \cr
		& &&&& I_{n-i-1} \cr
	} \text{ and  } \tau_i  \mapsto {}^ {\tilde{h}}  \mathcal{V}(\tau_i)=\bordermatrix {
		&  & i && i+1 &  \cr
		& I_{i-1} &  &&  &  \cr
		&  & 0 && \beta^{-1} &  \cr
		& & & & & \cr 
		& & \beta &&0 &  \cr
		&  &&&
		& I_{n-i-1} \cr}\text{ $ $ }
	\end{eqnarray}	
that does not factor through the welded braid group on $n$-strands. In contrast, ${}^{\tilde{h}}\phi_n$ factors through $\mathcal{WB}_n$ if and only if $\beta=\alpha$ and ${}^{\tilde{h}}\phi_n$ factors through $WB_n$ if and only if $\beta= 1.$	
\label{dernièreproposition}
\end{theo} 
\begin{proof}
It is a consequence of the relations of $VB_n$ to be verified by a direct calculation.
\end{proof}

\subsection{Generalized Gassner representation}
\label{sec:generalized Gassner}
Here, we propose an iterative procedure to find and generalize the extension of the Gassner representation of the pure Artin braid group on $n$ strands $P_n$. Recall the faithful representation $\psi:WB_n \rightarrow Aut(F_n)$ defined in section \ref{sec:two generalizations}. It follows from $WB_n=PW_n\rtimes S_n$ that $PW_n$ is embedded into the automorphism group $Aut(F_n)$ of the free group $F_n$. Then, using $PW_n \hookrightarrow Aut(F_n)$, a semidirect product $F_n\rtimes PW_n$ can be defined. Let $\overline{I}_{F_n}$ be the kernel of the canonical projection $\pi:\mathbb{C}[F_n\rtimes PW_n]\twoheadrightarrow \mathbb{C}[PW_n].$ The kernel $\overline{I}_{F_n}$ is called the relative augmentation ideal of $F_n$ in $\mathbb{C}[F_n\rtimes PW_n]$ and is generated by the elements $\{\xi_{n+1,1}-1,\xi_{n+1,2}-1,\dots,\xi_{n+1,n}-1\}$. For more details about the relative augmentation ideal, we refer to \cite{thetheoryofnilpotent,robinson}. Any generator $(\xi_{n+1,i}-1)$ of $\overline{I}_{F_n}$ can be written using the fundamental formula of the free differential calculus [\cite{fox},(2.3)]: 
\begin{eqnarray}
\xi_{n+1,i}-1= \sum \limits_{k=1} ^{n} \frac{\partial \xi_{n+1,i}}{ \partial \xi_{n+1,k}}( \xi_{n+1,k}-1) \text{  $$  $$  for all } i=1,\dots, n.
\end{eqnarray}
Excellent references for the Fox derivation are the Birman's book \cite{birman} and the Fox's original article \cite{fox}. Thanks to [Theorem 3.9, \cite{birman}], we define an action to the right of $PW_{n+1}$ on $\overline{I}_{F_n}$ by setting

\begin{eqnarray}
(\xi_{n+1,l}-1)\cdot\xi_{i,j}&=&\sum \limits_{k=1} ^{n} \frac{\partial (\xi_{n+1,l} \;\xi_{i,j})}{ \partial \xi_{n+1,k}}( \xi_{n+1,k}-1) \text{  $$  $$   } ( 1\leq l \leq n).
\label{Eq:actiondemccool}
\end{eqnarray}
By computing the coefficients $\left(\frac{\partial(\xi_{n + 1, l}\cdot \xi_{i, j})}{\partial\xi_{n + 1, k}} \right)_{l, k}$ of \eqref{Eq:actiondemccool}, we obtain
\begin{eqnarray}
(\xi_{n+1,l}-1)\cdot\xi_{i,j}=\sum \limits_{k=1} ^{n} \mathcal{C}(\xi_{i,j})_l^{k} \cdot( \xi_{n+1,k}-1)
\label{Eq:actiondemccoolenmatrice}
\end{eqnarray}
with 

\begin{eqnarray}
\mathcal{C}(\xi_{i,j})_l^{k}=\bordermatrix {
	&  & i  & &  & j &  \cr
	& \xi_{i,j}\cdot I_{i-1} &   & &  & \cr
	&  & \xi_{i,j}\cdot \xi_{n+1,j} &  & & \xi_{i,j} \cdot(1-\xi_{n+1,j}\cdot \xi_{n+1,i} \cdot\xi^{-1}_{n+1,j})& \cr
	&  & & &  & \cr
	& & 0 & & &\xi_{i,j}\cdot I_{j-i} &  \cr
	& & & & & & \cr
	&  & & &  & &\xi_{i,j}\cdot I_{n-j} \cr
} \text{  $ $   } (i<j)
\nonumber\\
\nonumber\\
\nonumber\\
\nonumber\\
\mathcal{C}(\xi_{i,j})_l^{k}=\bordermatrix {
	&  & j  & &  & i &  \cr
	& \xi_{i,j}\cdot I_{j-1} &   & &  & \cr
	&  & \xi_{i,j}\cdot I_{i-j}  &  & & 0& \cr
	&  & & &  & \cr
	& & \xi_{i,j} \cdot(1-\xi_{n+1,j}\cdot \xi_{n+1,i} \cdot\xi^{-1}_{n+1,j}) & & &\xi_{i,j}\cdot \xi_{n+1,j} &  \cr
	& & & & & & \cr
	&  & & &  & &\xi_{i,j}\cdot I_{n-i} \cr
} \text{  $ $   }(j<i) \nonumber\\
\label{Eq:matriceassociéeàmccool}
\end{eqnarray}
where $1\leq i \neq j \leq n.$ Moving from the generators $\xi_{i,j}$ in $PW_{n + 1}$ to their corresponding matrix $\mathcal{C}(\xi_{i,j})_l^{k}$, we find that the number of generators decreases by one, and then we assign each generator $\xi_{i,j}$ of $PW_n$ to this matrix. We arrive at the following result.
\begin{theo}
The mapping $\xi_{i,j} \mapsto \mathcal{C}(\xi_{i,j})$, $1\leq i\neq j\leq n$, determine a faithful matrix representation $\Psi_n:PW_n \rightarrow GL(n,\mathbb{C}[F_n\rtimes PW_n] )$ where
	\begin{eqnarray}
	\mathcal{C}(\xi_{i,j})&=&\bordermatrix {
		&  & i  & &  & j &  \cr
		& \xi_{i,j}\cdot I_{i-1} &   & &  & \cr
		&  & \xi_{i,j}\cdot \xi_{n+1,j} &  & & \xi_{i,j} \cdot(1-\xi_{n+1,j}\cdot \xi_{n+1,i} \cdot\xi^{-1}_{n+1,j})& \cr
		&  & & &  & \cr
		& &  0 & & &\xi_{i,j}\cdot I_{j-i} &  \cr
		& & & & & & \cr
		&  & & &  & &\xi_{i,j}\cdot I_{n-j} \cr
	}\text{  } (i<j)\nonumber\\
	\label{Eq:basisconj} 
	\end{eqnarray}
	and	
		
	\begin{eqnarray}
	\mathcal{C}(\xi_{i,j})&=&\bordermatrix {
		&  & j  & &  & i &  \cr
		& \xi_{i,j}\cdot I_{j-1} &   & &  & \cr
		&  & \xi_{i,j}\cdot I_{i-j}  &  & & 0& \cr
		&  & & &  & \cr
		& & \xi_{i,j} \cdot(1-\xi_{n+1,j}\cdot \xi_{n+1,i} \cdot\xi^{-1}_{n+1,j}) & & &\xi_{i,j}\cdot \xi_{n+1,j} &  \cr
		& & & & & & \cr
		&  & & &  & &\xi_{i,j}\cdot I_{n-i} \cr
	} \text{  $ $   }(j<i) \nonumber\\
	\label{Eq:basisconj1}
	\end{eqnarray}
\label{4.3.6}
\end{theo}
\begin{proof}
Since the faithful action $PW_n$ on $F_n$ extends over the free group algebra $\mathbb{C}[F_n]$ and, moreover, the matrix representation given in this theorem is equivalent to the right action of $PW_n$ on the augmentation ideal $I_{F_n}$ which is the kernel of the augmentation map $\epsilon:\mathbb{C}[F_n]\rightarrow \mathbb{C}$.
\end{proof}
\begin{remark}
Under Theorem \ref{4.3.6} and Formula \eqref{Eq:actiondemccoolenmatrice}, we can say that $\overline{I}_{F_n}$ satisfies the following abstract notation:
	\begin{eqnarray}
	\overline{I}_{F_n}\cdot \mathfrak{s}=\sum \mathcal{C}(\mathfrak{s})\cdot \overline{I}_{F_n}
	\label{Eq:écritureabstraite}
	\end{eqnarray}
	where $\mathfrak{s}$ represents the generators of $PW_{n+1}$ and $\mathcal{C}(\mathfrak{s})$ is their corresponding matrix considered as generators of $PW_n.$
\end{remark}

\begin{prop}
The homomorphism $\mathfrak{a}:\mathbb{C}[F_n\rtimes PW_n]\rightarrow \mathbb{C}[t^{\pm1}_1,\dots,t^{\pm1}_n]$ such that $\xi_{n+1,k}\mapsto t_k$ and $\xi_{i,j}\mapsto 1$ for all $1\leq k \leq n$ and $1\leq i \neq j \leq n.$ applied to each entry of the above matrix yields the representation ${}^{\mathfrak{a}}\Psi_n:PW_n\rightarrow Gl(n,\mathbb{C}[t^{\pm1}_1,\dots,t^{\pm1}_n]$, such that $\xi_{i,j}  \mapsto  {}^ {\mathfrak{a}}  \mathcal{C}(\xi_{i,j})$, where
\small{	\begin{eqnarray}
	{}^ {\mathfrak{a}}  \mathcal{C}(\xi_{i,j})=\bordermatrix {
		&  & i & & j &  \cr
		& I_{i-1} &  & &  & \cr
		&  & t_j &  &1-t_i&  \cr
		&  & &  & \cr
		& &  0 & & I_{j-i} &  \cr
		&  & & & & I_{n-j} \cr
	} (i<j)\text{ $ $ and   } {}^ {\mathfrak{a}}  \mathcal{C}(\xi_{i,j})=\bordermatrix {
	&  & j  & &  & i &  \cr
	& I_{j-1} &   & &  & \cr
	&  &  I_{i-j}  &  & &0& \cr
	&  & & &  & \cr
	& & 1-t_i & & &t_j &  \cr
	& & & & & & \cr
	&  & & &  & & I_{n-i} \cr
} (j<i)
	\end{eqnarray} }
with $1\leq i\neq j \leq n.$
\label{extensiondegassner}	
\end{prop}

\begin{remark}
The representation ${}^{\mathfrak{a}}\Psi_n$ ($n\geq1$) in Proposition \ref{extensiondegassner} was constructed respectively: by Bardakov \cite{barda05} using the Magnus construction in the Birman sense defined in \cite{birman} and by Rubinsztein \cite{rubinsztein} using the construction in the original sense considered by Magnus in \cite{magnussurvey}. This representation is not faithful for $n\geq 2$ and is an extension of Gassner representation of the pure braid group $P_n$ (see \cite{barda05}). In addition, M. Nasser and M. Abdulrahim recently studied the irreducibility of ${}^{\mathfrak{a}}\Psi_n$ in \cite{nasserabdul}. 
\end{remark}
Let us now see how we can iterate the above procedure as follows. We first iteratively construct semi-direct products
\begin{equation}
\begin{split}
PW_n(1) &:= F_{n} \rtimes PW_n \leq PW_{n+1}\\
PW_n(2) &:= F_{n} \rtimes( F_{n-1} \rtimes PW_{n-1}) \leq PW_{n+1}\\
\vdots\\
PW_n(r) &:= F_{n} \rtimes(F_{n-1} \rtimes( \dots ( F_{n+1-r} \rtimes PW{n+1-r})\dots )) \leq PW_{n+1} \text{  $ $   } (1\leq r\leq n),
\end{split}
\label{Eq:4.54}
\end{equation}
where $F_k$ is a free group on $\{\xi_{k+1,1},\dots,\xi_{k+1,k}\}$ (see \cite{barda03,cpvw}). From the decompositions \eqref{Eq:4.54}, we iteratively define relative augmentation ideals $\overline{I}_{F_{n+1-r}}$ associated with the free groups $F_{n + 1-r}$, and then consider the action of $PW_{n + 1}$ on the direct product $\mathfrak{d}_{l^r}\cdot\mathfrak{d}_{l^{r-1}}\cdots\mathfrak{d}_{l^2}\cdot\mathfrak{d}_{l ^ 1}$, where $\mathfrak{d}_{l^s}:=(\xi_{n-r+2,l^s}-1)$ for each $s=1,\dots,r$. Here the $\xi_{i,j}$-action ($1\leq i \neq j \leq n + 1$) on the direct product $\mathfrak{d}_{l^r}\cdot\mathfrak{d}_{l^{r-1}}\cdots\mathfrak{d}_{l^2}\cdot\mathfrak{d}_{l ^ 1}$ is given by iteration of \eqref{Eq:actiondemccool} as follows:
\footnotesize{\begin{eqnarray}	\mathfrak{d}_{l^r}\cdot\mathfrak{d}_{l^{r-1}}\cdots\mathfrak{d}_{l^2}\cdot \mathfrak{d}_{l^1}\cdot \xi_{i,j}&=&\mathfrak{d}_{l^r}\cdot\mathfrak{d}_{l^{r-1}}\cdots\mathfrak{d}_{l^{2}}\cdot\sum \limits_{k^1=1} ^{n} \mathcal{C}(\xi_{i,j})_{l^1}^{k^1} \cdot\mathfrak{d}_{k^1}\text{  $ $ by  } \eqref{Eq:actiondemccoolenmatrice} \nonumber\\	&=& \mathfrak{d}_{l^r}\cdot\mathfrak{d}_{l^{r-1}}\cdots \mathfrak{d}_{l^3}\cdot\sum \limits_{k^{2}=1} ^{n-1}\sum \limits_{k^1=1} ^{n} \mathcal{C}(\mathcal{C}(\xi_{i,j})_{l^1}^{k^1})_{l^{2}}^{k^{2}}\cdot\mathfrak{d}_{k^1}\cdot\mathfrak{d}_{k^{2}}\text{  $ $ by } \eqref{Eq:écritureabstraite} \text{  $ $ and } \eqref{Eq:actiondemccoolenmatrice} \nonumber\\	&\vdots&\nonumber\\
	&=&\sum \limits_{k^r=1} ^{n+1-r}\cdots \sum \limits_{k^{2}=1} ^{n-1}\sum \limits_{k^1=1} ^{n} \mathcal{C}(\cdots\mathcal{C}(\mathcal{C}(\xi_{i,j})_{l^1}^{k^1})_{l^{2}}^{k^{2}}\cdots)_{l^r}^{k^r}\cdot \mathfrak{d}_{k^1}\cdots\mathfrak{d}_{k^r} \text{ by iteration of } \eqref{Eq:écritureabstraite} \text{ and } \eqref{Eq:actiondemccoolenmatrice},\nonumber\\
	\end{eqnarray} }
where the coefficients $\mathcal{C}(\mathcal{C}(\xi_{i,j})_{l^1}^{k^1})_{l^{2}}^{k^{2}}$ are exactly like the matrices \eqref{Eq:matriceassociéeàmccool}, but the matrix elements $\xi_{i,j}\in PW_{n+1}$ present in these matrices are replaced by their corresponding matrix representations $\mathcal{C}(\xi_{i,j})_{l^1}^{k^1}.$ Similarly, $\mathcal{C}(\cdots\mathcal{C}(\mathcal{C}(\xi_{i,j})_{l^1}^{k^1})_{l^{2}}^{k^{2}}\cdots)_{l^r}^{k^r}$ are also matrices obtained iteratively from the previously introduced matrices by iteratively replacing the matrix elements present in \eqref{Eq:matriceassociéeàmccool} with their corresponding matrix representations. The substitutions of $\xi_{i,j}$ by $\mathcal{C}(\xi_{i, j})_{l^1}^{k^1}$ reduce the number of generators of $PW_{n + 1}$ by two and, moreover, since these substitutions respect the McCool relations \eqref{Eq:2.7}, it must be clear that the mapping $\xi_{i,j}\mapsto \mathcal{C}(\mathcal{C}(\xi_{i,j})_{l^1}^{k^1})_{l^{2}}^{k^{2}}$ establishes a matrix representation of $PW_{n-1}$ in $GL(2(n-1), \mathbb{C}[PW_n(2)]).$ For the same reason as above, the mapping $\xi_{i,j}\mapsto \mathcal{C}(\cdots\mathcal{C}(\mathcal{C}(\xi_{i,j})_{l^1}^{k^1})_{l^{2}}^{k^{2}}\cdots)_{l^r}^{k^r}$ defines a representation of $PW_{n+1-r}$ in $GL(r(n+1-r),\mathbb{C}[PW_n(r)]).$ To illustrate more precisely the form of this matrix representation, if we use $\xi^{(r)}_{i, j}$ to represent the matrix representation of $\xi_{i,j}$ at the $r$-th iteration, then the matrix representation $\Psi^{(r)}_n$ of $PW_{n + 1-r}$ in  $GL(r(n+1-r),\mathbb{C}[PW_n(r)])$ is given by
\footnotesize{ \begin{eqnarray}
	\xi^{(0)}_{i,j}&=&\xi_{i,j} \in PW_{n+1} \text{ \color{blue}  $$ $$ (\emph{ initial representation } ) }
	\nonumber\\
	\xi^{(r)}_{i,j}&=&\bordermatrix {
		&  & i  & &  & j &  \cr
		& \xi^{(r-1)}_{i,j}\cdot I_{i-1} &   & &  & \cr
		&  & \xi^{(r-1)}_{i,j}\cdot t^{(r)}_j \cdot \xi^{(r-1)}_{n+1,j} &  & & \xi^{(r-1)}_{i,j} \cdot(I^{(r)}-t^{(r)}_i\cdot\xi^{(r-1)}_{n+1,j}\cdot \xi^{(r-1)}_{n+1,i} \cdot\xi^{-(r-1)}_{n+1,j})& \cr
		&  & & &  & \cr
		&  & & &  & \cr
		& & 0 & & &\xi^{(r-1)}_{i,j}\cdot I_{j-i} &  \cr
		& & & & & & \cr
		&  & & &  & &\xi^{(r-1)}_{i,j}\cdot I_{n+1-r-j}\cr
	}(i<j)\nonumber\\
	\label{Eq:itérationdemccool}
	\end{eqnarray} }
and
\footnotesize{\begin{eqnarray}
	\xi^{(r)}_{i,j}&=& \bordermatrix {
		&  & j  & &  & i &  \cr
		&\small{ \xi^{(r-1)}_{i,j}\cdot I_{j-1} } &   & &  & \cr
		&  & \xi^{(r-1)}_{i,j}\cdot I_{i-j}  &  & & 0& \cr
		&  & & &  & \cr
		& &\small{ \xi^{(r-1)}_{i,j} \cdot(I^{(r)}-t^{(r)}_i\cdot\xi^{(r-1)}_{n+1,j}\cdot \xi^{(r-1)}_{n+1,i} \cdot\xi^{-(r-1)}_{n+1,j}) }& & &\small {\xi^{(r-1)}_{i,j}\cdot t^{(r)}_j \cdot \xi^{(r-1)}_{n+1,j} } &  \cr
		& & & & & & \cr
		&  & & &  & &\small {\xi^{(r-1)}_{i,j}\cdot I_{n+1-r-i} } \cr
	} \text{   }(j<i)\nonumber\\
	\label{Eq:itérationdemccool1}
	\end{eqnarray} }
where $r=1,\dots,n$, $t^{(r)}_j,t^{(r)}_i$ are nonzero, not necessarily distinct, parameters in $\mathbb{C}$ introduced at each iteration of each matrix representation, and $I^{(r)}$ is the identity matrix of the size of the recurrent matrices $\xi^{(r-1)}_{i,j}$. Applying them to a trivial representation $\xi^{(0)}_{i,j}=1$ for all $1\leq i\neq j\leq n+1$ of $PW_{n+1},$ the formulas \eqref{Eq:itérationdemccool}-\eqref{Eq:itérationdemccool1} yield for the first iteration, the extension of the Gassner representation ${}^ {\mathfrak{a}} \Psi_n$ of $P\Sigma_n$ (see Proposition \ref{extensiondegassner}) and it yields in the second iteration, a product tensor ${}^ {\mathfrak{a}} \Psi_n\otimes {}^ {\mathfrak{a}} \Psi_n$ of two extensions ${}^ {\mathfrak{a}} \Psi_n$ of the Gassner representation. Let us illustrate this with a simple example.
\begin{example}
Applying them to a trivial representation $\xi^{(0)}_{i,j}=1$ for all $1\leq i\neq j\leq 5$ of $P\Sigma_{5},$ the formulas \eqref{Eq:itérationdemccool}-\eqref{Eq:itérationdemccool1} yield, without any real surprise, for the first iteration the representation ${}^ {\mathfrak{a}} \Psi_4$ of $P\Sigma_4.$ For the second iteration $r=2,$ we obtain for example by applying the formulas \eqref{Eq:itérationdemccool}-\eqref{Eq:itérationdemccool1} on ${}^ {\mathfrak{a}} \Psi_4(\xi_{1,2})$:
	\begin{eqnarray}
	\xi^{(2)}_{1,2}= \bordermatrix{ &  & & & & & & & & & & & \cr
		& s_2t_2 & s_2(1-t_1) & 0 & 0 & t_2(1-s_1) & (1-t_1)(1-s_1) & 0 & 0 & 0 & 0 & 0 &0\cr
		& 0 & s_2 & 0 & 0 & 0 & 1-s_1 & 0 & 0 & 0 & 0 & 0 &0 \cr
		& 0 & 0 & s_2 & 0 & 0 & 0 & 1-s_1 & 0 & 0 & 0 & 0  &0\cr
		& 0 & 0 & 0 &s_2t_2 & 0 & 0 & 0 & 1-s_1t_1 & 0 & 0 & 0 & 0 \cr
		& 0 & 0 & 0 & 0 & t_2 & 1-t_1 & 0 & 0 & 0 & 0 & 0 & 0 \cr
		& 0 & 0 & 0 & 0 & 0 & 1 & 0 & 0 & 0 & 0 & 0 & 0 \cr
		& 0 & 0 & 0 & 0 & 0 & 0 & 1 & 0 & 0 & 0 & 0 & 0 \cr
		& 0 & 0 & 0 & 0 & 0 & 0 & 0 & 1 & 0 & 0 & 0 & 0 \cr
		& 0 & 0 & 0 & 0 & 0 & 0 & 0 & 0 & t_2 & 1-t_1 & 0 & 0 \cr
		& 0 & 0 & 0 & 0 & 0 & 0 & 0 & 0 & 0 & 1 & 0 & 0 \cr
		& 0 & 0 & 0 & 0 & 0 & 0 & 0 & 0 & 0 & 0 & 1 & 0 \cr
		& 0 & 0 & 0 & 0 & 0 & 0 & 0 & 0 & 0 & 0 & 0 & 1 \cr
	}\text{  }
	\end{eqnarray}
	As we can see by omitting columns 4,8 and 12 and their rows in the above matrix, the result of $\xi^{(2)}_{1,2}$ becomes a tensor product of two matrices $\xi^{(1)}_{1,2}:$
	
	\begin{eqnarray}
	\widehat{\xi^{(2)}}_{1,2}= \bordermatrix{ &  & & & & & & & & & & & \cr
		& s_2t_2 & s_2(1-t_1) & 0 & t_2(1-s_1) & (1-t_1)(1-s_1) & 0 & 0 & 0 & 0 \cr
		& 0 & s_2 & 0  & 0 & 1-s_1 & 0  & 0 & 0 & 0 \cr
		& 0 & 0 & s_2  & 0 & 0 & 1-s_1 & 0 & 0 & 0 \cr
		& 0 & 0 & 0  & t_2 & 1-t_1 & 0 & 0 & 0 & 0 \cr
		& 0 & 0 & 0 & 0 & 1 & 0 & 0 & 0 & 0 \cr
		& 0 & 0 & 0 & 0 & 0 & 1  & 0 & 0 & 0 \cr
		& 0 & 0 & 0 & 0 & 0 & 0  & t_2 & 1-t_1 & 0  \cr
		& 0 & 0 & 0 & 0 & 0 & 0  & 0 & 1 & 0  \cr
		& 0 & 0 & 0 & 0 & 0 & 0  & 0 & 0 & 1  \cr
	}=\xi^{(1)}_{1,2}\otimes \xi^{(1)}_{1,2}. \nonumber\\
	\end{eqnarray}
	In other words, it is a tensor product $\left({}^{\mathfrak{a}}\Psi_3\otimes {}^{\mathfrak{a}}\Psi_3 \right)(\xi_{1 ,2})$ of $P\Sigma_3.$ 
\end{example}

\bibliographystyle{abbrv}	

\end{document}